\documentclass[a4paper]{amsart}

\usepackage[mathscr]{eucal}
\usepackage{amsmath}
\usepackage{amsfonts}
\usepackage{amssymb}
\usepackage{footnpag}
\usepackage[dvips]{graphicx}

\theoremstyle{plain}

\newtheorem{thm}{Theorem}[section]
\newtheorem{prop}{Proposition}[section]

\newtheorem{rem}{Remark}[section]

\newtheorem{cor}{Corollary}[section]

\newcommand{\C}{\mathbb{C}}

\newcommand{\F}{\mathbb{F}}

\begin{document}

\title{On the nonexistence of ternary extremal self-dual codes}
\author{Tsuyoshi Miezaki}

\maketitle \vspace{-0.2in}
\begin{center}
Graduate School of Mathematics, Kyushu University\\
Hakozaki 6-10-1 Higashi-ku, Fukuoka, 812-8581 Japan\\ \quad
\end{center} \vspace{0.1in}

\begin{quote}
{\small\bfseries Abstract.}

In this note, we give a new nonexistence result of ternary extremal 
self-dual codes. 

\vspace{-0.15in}

\quad
\end{quote}

\section{Introduction}                                   
Let $\F_{q}$ be a finite field of order $q$ and $q$ is a prime or prime power. A linear code $C$ of length $n$ is a linear subspace of $\F_{q}^{n}$. We define the inner product by $\langle \bold{x}|\bold{y}\rangle = x_1 y_1 + \cdots + x_n y_n$, where $\bold{x} =(x_1, \ldots, x_n)$ and $\bold{y} = (y_1, \ldots, y_n)$. Then, the dual $C^{\perp}$ of a linear code $C$ is defined as follows: $C^{\perp}=\{ \bold{y}\in \F_{q}^{n}\ | \ \langle \bold{x}|\bold{y}\rangle =0 \text{ for all } \bold{x}\in C\}$. A linear code $C$ is called self-dual if and only if $C=C^{\perp}$. The weight $wt(\bold{x})$ is the number of its nonzero components. The weight enumerator of a code $C$ is
\begin{align*}
W(x, y)&=\sum_{\bold{u}\in C}x^{n-wt(\bold{u})}y^{wt(\bold{u})}=x^{n}+\sum_{i=1}^{n}A_{i}x^{n-i}y^{i}, 
\end{align*}
where $wt(\bold{u})$ denotes the weight of $\bold{u}$ and $A_{i}$ is the number of codewords of weight $i$. 

A code over $\F_{3}$ is called ternaty. In this note, we consider ternary self-dual codes. A ternary self-dual code exists if and only if the length $n$ is a multiple of $4$. 
It is well known that the weight enumerator $W(x, y)$ is in the space $\C[f, g]$, where $f=x^4+8xy^3$, $g=y^3(x^3-x^3)^3$, \cite{1}. 
Then, the minimum weight $d$ of a ternary self-dual code of length $n$ satisfies the Mallows--Sloane bound, $d\leq 3[n/12]+3$, \cite{1}. A self-dual code $C$ is called extremal if the minimum weight attain the Mallows--Sloane bound. 

We set $n=4j=12\alpha +4\nu$, $0\leq \nu\leq 2$. Then, we can choose the $a_{1}, \ldots, a_{\alpha}$ so that 
\begin{align*}
W^{\ast}(x, y)=\sum_{r=0}^{\alpha}a_{r}f^{j-3r}g^{r}=x^{n}+\sum_{r=\alpha +1}^{[n/3]}A^{\ast}_{3r}x^{n-3r}y^{3r}. 
\end{align*}
This is called an extremal weight enumerator and a code having this weight enumerator is a ternary extremal code of length $n$. If there exists a negative coefficient of $W^{\ast}(x, y)$, then such an extremal code does not exist. In fact, the following results hold. 

\begin{thm}[cf. \cite{1}]\label{thm:MS}
For extremal ternary self-dual codes of length $n$, the coefficient of the highest power of $y$ is negative for $n=24i$ $($$i\geq 3$$)$, $24i+4$ $($$i\geq 7$$)$, and the next-to-highest coefficient is negative for $n=24i+12$ $($$i\geq11$$)$. 
\end{thm}

\begin{thm}[cf. \cite{3}]\label{thm:Z}
For extremal ternary self-dual codes, $A^{\ast}_{3(\alpha +2)}$ $($the third nonzero coefficient in $W^{\ast}(x, y)$$)$ is negative if and only if $n=12i$ $($$i\geq 70$$)$, $n=12i+4$ $($$i\geq 75$$)$, $n=12i+8$ $($$i\geq 78$$)$. 
\end{thm}

In \cite{3}, Zhang remarked that there are some missing cases in Theorem \ref{thm:MS}, which are $n=24i+8$, $24i+16$ and $24i+20$. Moreover, by Theorem \ref{thm:Z}, ternary extremal self-dual codes of length $n$ greater than $944$ do not exist. In this note, we consider these three cases and give the length $n$ up to $944$ for which weight enumerator has some negative coefficient. In other words, we show the nonexistence of ternary self-dual codes for the length $n$ given later. 

\section{Nonexistence of extremal self-dual codes}
In \cite{2}, Mallows and Sloane gave the explicit form of the extremal weight enumerator of ternary self-dual codes. 

\begin{thm}[cf. \cite{2}]\label{thm:MS2}
We set $n=4j=12\alpha +4\nu$, $0\leq \nu\leq 2$. The extremal weight enumerator of a ternary self-dual codes of length $n$ is given by 
\begin{align*}
W^{\ast}(x, y)=\sum_{r=0}^{\alpha}a_{r}f^{j-3r}g^{r}, 
\end{align*}
where $a_{0}=1$ and $a_{r}$, $1\leq r\leq \alpha$, is equal to 
\begin{align*}
\frac{j}{r}\sum_{i=0}^{r-1}(-8)^{i+1} \binom{j-3r+i}{i} \binom{4r-i-2}{r-i-1}. 
\end{align*}
\end{thm}

Then, the coefficient of the highest power of $y$ and the next-to-highest coefficient are computable explicitly using Theorem \ref{thm:MS2}. 

\begin{cor}\label{cor:main}
Let the notation be the same as Theorem \ref{thm:MS2} and we set 
\begin{align*}
W^{\ast}(x, y)=x^{n}+\sum_{r=1}^{[n/3]}A^{\ast}_{3r}, 
\end{align*}
where $A^{\ast}_{3}=\cdots =A^{\ast}_{\alpha }=0$. 
Then 
\begin{align*}
\left\{
\begin{array}{lll}
A^{\ast}_{[n/3]}&=64 a_{\alpha}\ &{\rm if}\ n\equiv 8\pmod {24}\\
A^{\ast}_{[n/3]-1}&=4096a_{\alpha-1}+(2n-9)a_{\alpha}
\ &{\rm if}\ n\equiv 16\pmod {24} \\
A^{\ast}_{[n/3]}&=-64 a_{\alpha}\ &{\rm if}\ n\equiv 20\pmod {24}. 
\end{array}
\right. 
\end{align*}
\end{cor} 

\begin{proof}
The coefficient of the highest power of $y$ depends on the term 
$a_{\alpha}f^{j-3\alpha}g^{\alpha}$ and the next-to-highest coefficient 
depends on the terms $a_{\alpha-1}f^{j-3(\alpha-1)}g^{\alpha-1}+a_{\alpha}f^{j-3\alpha}g^{\alpha}$. Then the results are easy calculation. 
\end{proof}

We calculated the extremal weight enumerator of length up to 944 
using Corollary \ref{cor:main}. 
Then, we verified the following fact. 
\begin{prop}\label{prop:mainprop}
For extremal ternary self-dual codes of length $n$, the coefficient of the highest power of $y$ is negative for $n=24i+8$ $($$11\leq i\leq 38$$)$, $24i+20$ $($$19\leq i\leq 38$$)$, and the next-to-highest coefficient is negative for $n=24i+16$ $($$15\leq i\leq 36$$)$. 
\end{prop}

In Appendix, for each case, we give details for the minimum length for which weight enumerator has the negative coefficient, that is, $n=272$, $376$, $476$. 

By Proposition \ref{prop:mainprop} together with Theorems \ref{thm:MS}, \ref{thm:Z}, we get the following bounds for the length $n$ of ternary self-dual codes. 
\begin{thm}\label{thm:main}
For extremal ternary self-dual codes of length $n$, the coefficient of the highest power of $y$ is negative if and only if $n=24i$ $($$i\geq 3$$)$, $24i+4$ $($$i\geq 7$$)$, $n=24i+8$ $($$i\geq 11$$)$, $24i+20$ $($$i\geq 19$$)$, and the next-to-highest coefficient is negative if and only if $n=24i+12$ $($$i\geq11$$)$, $n=24i+16$ $($$i\geq 15$$)$. So, an extremal ternary self-dual code does not exist for $n=24i$ $($$i\geq 3$$)$, $24i+4$ $($$i\geq 7$$)$, $n=24i+8$ $($$i\geq 11$$)$, $n=24i+12$ $($$i\geq11$$)$, $n=24i+16$ $($$i\geq 15$$)$, $24i+20$ $($$i\geq 19$$)$. 
\end{thm}

\begin{rem}
We remark that for ternary seld-dual codes, all the coefficients of the extremal weight enumerator are positive for the length $n$ which are not mentioned in Theorem \ref{thm:main}. We verify this fact by a computer calculation $($using Mathematica$)$. 
\end{rem}

\begin{rem}
In \cite{3}, Zhang shows that the third coefficient of the extremal weight enumerator of a doubly even self-dual codes of length $n$ is negative for $n=24i$ $($$i\geq 154$$)$, $n=24i+8$ $($$i\geq 159$$)$, $n=24i+16$ $($$i\geq 164$$)$. We remark that all the coefficients of the extremal weight enumerator are positive for $n=24i$ $($$i\leq 153$$)$, $n=24i+8$ $($$i\leq 158$$)$, $n=24i+16$ $($$i\leq 163$$)$. We verify this fact by a computer calculation $($using Mathematica$)$. 
\end{rem}

\begin{center}
{\bf Acknowledgment}
\end{center}
The author is supported by JSPS Research Fellowship. The author would like to thank Eiichi Bannai, Masaaki Harada, Akihiro Munemasa for their helpful discussions on this research. 

\begin{appendix}
\section{$n=272$, $376$, $476$}
\hspace{-12pt}(i)\ $n\equiv 8\pmod {24}$
\begin{align*}
\begin{array}{l|l}
n& a_{\alpha}, 
A^{\ast}_{[n/3]} \\ \hline 
272 & a_{22}= -7872953986244054304, \\
&A^{\ast}_{270}=-503869055119619475456\\
296 & a_{24}=-2744029932081068834640, \\
&A^{\ast}_{294}=-175617915653188405416960\\
320 & a_{26}=-315678406695936752816640, \\
&A^{\ast}_{318}=-20203418028539952180264960\\
344 & a_{28}=-30299903951133235484744224, \\
&A^{\ast}_{342}=-1939193852872527071023630336\\ 
368 & a_{30}=-2741343458082960767366977664, \\
&A^{\ast}_{366}=-175445981317309489111486570496 \\
392 & a_{32}=-242684023678935047739004262256, \\
&A^{\ast}_{390}=-15531777515451843055296272784384 \\
416 & a_{34}=-21317169142068615385452966769920, \\
&A^{\ast}_{414}=-1364298825092391384668989873274880 \\
440 & a_{36}=-1868185741186233164549296424639680, \\
&A^{\ast}_{438}=-119563887435918922531154971176939520 \\
464 & a_{38}=-163702510038923775966285120477012096, \\
&A^{\ast}_{462}=-10476960642491121661842247710528774144 \\
488 &a_{40}=-14354535735664486040811516311170482368, \\
&A^{\ast}_{486}=-918690287082527106611937043914910871552 \\
512 &a_{42}=-1259902182653364764424148281390773551104, \\
&A^{\ast}_{510}=-80633739689815344923145490009009507270656 \\
536 &a_{44}=-110692831895932549098687986745256783062960, \\
&A^{\ast}_{534}=-7084341241339683142316031151696434116029440 \\
560 &a_{46}=-9734710458637992092026263038211311974559040, \\
&A^{\ast}_{558}=-623021469352831493889680834445523966371778560 \\
\end{array}
\end{align*}

\begin{align*}
\begin{array}{l|l}
n& a_{\alpha}, 
A^{\ast}_{[n/3]} \\ \hline 
584 &a_{48}=-856883483582797698500319427732292467760832912, \\
&A^{\ast}_{582}=-54840542949299052704020443374866717936693306368 \\
608 &a_{50}=-75489699451488581447144688357236834925309686016, \\
&A^{\ast}_{606}=-4831340764895269212617260054863157435219819905024 \\
632 &a_{52}=-6655702067690834059974283743113198460090480573120, \\
&A^{\ast}_{630}=-425964932332213379838354159559244701445790756679680 \\
656 &a_{54}=-587239823791971638229294202183655622522374288738240, \\
&A^{\ast}_{654}=-37583348722686184846674828939753959841431954479247360 \\
680 &a_{56}=-51847772061205519283387706376624194411645827116923360, \\
&A^{\ast}_{678}=-3318257411917153234136813208103948442345332935483095040 \\
704 &a_{58}=-4580547570370715126128873365093950442946918553754776576, \\
&A^{\ast}_{702}=-293155044503725768072247895366012828348602787440305700864 \\
728 &a_{60}=-404910564817501337290382547302278935209903073669968549696, \\
&A^{\ast}_{726}=-25914276148320085586584483027345851853433796714877987180544 \\
752 &a_{62}=-35812835197972512952024947381110335224799769173725029839104, \\
&A^{\ast}_{750}=-2292021452670240828929596632391061454387185227118401909702656\\
776 &a_{64}=-3169138428633901231181913977283450477788995635052821390015600, \\
&A^{\ast}_{774}=-202824859432569678795642494546140830578495720643380568960998400\\
800 &a_{66}=-280577613193127107157764299638996491155166379020824310067308800, \\
&A^{\ast}_{798}=-17956967244360134858096915176895775433930648257332755844307763200 \\
824 &a_{68}=-2485202165883707951985598734777598604314309255294697585685262636
8, \\
&A^{\ast}_{822}=-
1590529386165573089270783190257663106761157923388606454838568087552 \\
848 &a_{70}=-2202196148383627422274663559827948820347657470848878126067954281
344, \\
&A^{\ast}_{846}=-
140940553496552155025578467828988724502250078134328200068349074006016 \\
872 &a_{72}=-1952205178849727649791721360107420687947671575339491214112164662
53888, \\
&A^{\ast}_{870}=-
1249411314463825695866701670468749240286509808217274377031785384024883
2 \\
896 &a_{74}=-1731253139633577068357222419311556274623038029710586077404740969
8222080, \\
&A^{\ast}_{894}=-
1108002009365489323748622348359396015758744339014775089539034220686213
120 \\
920 &a_{76}=-1535863295223138298753646585853346282622456158809348646029139051
551145920, \\
&A^{\ast}_{918}=-9829525089428085112023338149461416208783719416379831334586489\
9299273338880 \\
944 &a_{78}=-1362992056754976406905872968762068822037202029564428321228680204
62858775808, \\
&A^{\ast}_{942}=-8723149163231849004197587000077240461038092989212341255863553
309622961651712 
\end{array}
\end{align*}

\hspace{-12pt}(ii)\ $n\equiv 16\pmod {24}$
\begin{align*}
\begin{array}{l|l}
n& a_{\alpha-1}, a_{\alpha}, 
A^{\ast}_{[n/3]-1} \\ \hline 
376 &a_{30}=11261493712193575449464340416, \\
&a_{
31}=-79801722496094302883416289152, \\
&A^{\ast}_{372}=-13165601569453182001372364496000\\
400 &a_{32}=441888898139839668601207079200, \\
&a_{
33}=-6950256276649549829449841225600, \\
&A^{\ast}_{396}=-3687675788049010632504280213046400\\
424 &a_{34}=16394113470331217869141611757248, \\
&a_{
35}=-606519885441368468657517572782656, \\
&A^{\ast}_{420}=-441719895110831476811653201806960576\\
448 &a_{36}=536274881987240455245925701691904, \\
&a_{
37}=-53018454175840507150162030554867712, \\
&A^{\ast}_{444}=-44830786937350792937506409428037621760\\
\end{array}
\end{align*}

\pagebreak
\vspace{10pt}(ii)\ $n\equiv 16\pmod {24}$
\vspace{-11pt}
\begin{align*}
\begin{array}{l|l}
n& a_{\alpha-1}, a_{\alpha}, 
A^{\ast}_{[n/3]-1} \\ \hline 
472 &a_{38}=11186655811029329434543269718219840, \\
&a_{39}=-4641517355814653901038572806560252288, \\
&A^{\ast}_{468}=-4293998185484725264107176341368007424640\\
496 &a_{40}=-427807284262010340700050571888614528, \\
&a_{41}=-406888023028715569042794027925684406784, \\
&A^{\ast}_{492}=-401723225273564598724573936593403536975360\\
520 &a_{42}=-92446581307040731231409486276307429120, \\
&a_{43}=-35712049812783023444078610813331049405040, \\
&A^{\ast}_{516}=-37197784554012936005968901004332067166271760\\
544 &a_{44}=-10245925750546386007769026911283939159104, \\
&a_{465}=-3137842254787954331766805187499724994341632, \\
&A^{\ast}_{540}=-3427699104790440721064204731540822283690310912\\
568 &a_{46}=-982947450888674233355426028486625484972832, \\
&a_{47}=-275982821766429030417241211253764949553992384, \\
&A^{\ast}_{564}=-315058792889605526940054670095674316133798136640\\
592 &a_{48}=-89693252321825420600496997362770544679251360, \\
&a_{49}=-24295796231535727604680318622890951981427680896, \\
&A^{\ast}_{588}=-28914944133564676858279010083094776729183738623360\\
616 &a_{50}=-8026912927711453941452521598868703485590141376, \\
&a_{51}=-2140653278232020640128117056850678924273791411008, \\
&A^{\ast}_{612}=-2650897194629667358220876688997346533863824114738880\\
640 &a_{
52}-712877944212068232540658862265006546757418664960, \\
&a_{
53}-188755374623690299975151023348103049880929304965120, \\
&A^{\ast}_{636}=-242828029206203002748903489375276443214179533462343680\\
664 &a_{
54}=-63131244080866258996313898504641769228020507849120, \\
&a_{
55}=-16655780326477032641903955663011085918793336423214912,\\ 
&A^{\ast}_{660}=-22227559826378434251520219247786635013646382742370464448\\
688 &a_{
56}=-5585948312654604444473291096286362330072408852296256, \\
&a_{
57}=-1470689084857734214358876519897473520883792783898787072,\\ 
&A^{\ast}_{684}=-2033312023289155930833146803030235243152121322248647392000\\
712 &a_{
58}=-494225270561367320028724723939101928067625312541795200, \\
&a_{
59}=-129941201445319285129486987492297441687572118692952646976,\\ 
&A^{\ast}_{708}=-185891146753346149001061743770855441485279541230699188610240\\
736 &a_{
60}=-43738879452320934005955894801832206861158358165595289856, \\
&a_{
61}=-11487484164029988549910683533790938795012660764252345740288,\\ 
&A^{\ast}_{732}=-1698534378221257979420772535504444817640682733314746012529152
0\\
760 &a_{
62}=-3872358905211529931725585546555728885910410725168282924160, \\
&a_{
63}=-1016105224553602100116684019450350373484875195904456495947520,\\ 
&A^{\ast}_{756}=-1551196176376239199876657551788171679852335463341923052234062
080\\
784 &a_{
64}=-342975399724365658008866745753091264330623560747127210011616, \\
&a_{
65}=-89923533208716771622371860893946418844699213766843605684714368, \\
&A^{\ast}_{780}=-1415956155096604486944820493242671287975843083673294143146772\
78848\\
808 &a_{
66}=-30389906591723821070046717068350880583117383557269660184749888, \\
&a_{
67}=-7961877992160201363477815938406804924854697772744893544139795776,\\
 &A^{\ast}_{804}=-
12919214990801144362211761566131700721109948123851620453549387353280\\
832 &a_{
68}=-2693821891332165338133436316145155487452402042147479219740572160,\\
 &a_{
69}=-7052659222061828538508726417870958604612412679358082870582131978
24,\\ &A^{\ast}_{828}=-
1178248995718129172348188777308574205939959337198398789965400225966080\\
\end{array}
\end{align*}

\pagebreak
\vspace{10pt}(ii)\ $n\equiv 16\pmod {24}$
\vspace{-11pt}
\begin{align*}
\begin{array}{l|l}
n& a_{\alpha-1}, a_{\alpha}, 
A^{\ast}_{[n/3]-1} \\ \hline 
856 &a_{
70}=-2388766983717336918141061720759236479204556018229319202355162880
64,\\ &a_{
71}=-6249917852664063370567737485453048899893833432747230055129898142
5792,\\ &A^{\ast}_{852}=-
1074145399873996204024391482580884060270741695047520569841468400840339
20\\
880 &a_{
72}=-2119024806439757178660594845900947475615108617650824410368719030
2080, \\&a_{
73}=-5540761864836300800626631791799135532179367176071433248381459697
180160, \\&A^{\ast}_{876}=-
9788669281400135155935170232328389125447266774280057385764638661239779
840\\
904 &a_{
74}=-1880389494832103699293661456701439577327212947871414981826519388
882176, \\&a_{
75}=-4913930400162537256399704675842273828431908296188004964978039826
28115648, \\&A^{\ast}_{900}=-
89171815436007274917861370851067415824363256671870340896511078816484]14
43648
\end{array}
\end{align*}

\hspace{-12pt}(iii)\ $n\equiv 20\pmod {24}$
\begin{align*}
\begin{array}{l|l}
n& a_{\alpha}, 
A^{\ast}_{[n/3]} \\ \hline 
476 &a_{39}=-1534201729056210732510378215650788800, \\
&A^{\ast}_{471}=-1477107137796128531650000385882276746240\\
500 &a_{41}=-134515010156247677274376431457788336000, \\
&A^{\ast}_{495}=-661923588246447626557238419517002415155200\\
524 &a_{43}=-11809682860520431518918202868161515902568, \\
&A^{\ast}_{519}=-82318278638696194278518292303059769395852160\\
548 &a_{45}=-1037995066725346452142518957162332543754336, \\
&A^{\ast}_{543}=-8428238763638643560243884563137644676299792896\\
572 &a_{47}=-91323929307918317042299735831678713063480672, \\
&A^{\ast}_{567}=-808309583260181473530593364281967075843110969856\\
596 &a_{49}=-8042005937180513916472095265588891227006959520, \\
&A^{\ast}_{591}=-75521483119365401801437288380702086144565783707136 \\
620 &a_{
51}-708763495653366966463564758903522391980006841440, \\
&A^{\ast}_{615}=-6977184797681093942233761765095086754438283587719680\\
644 &a_{
53}=-62512560340137400775498728136754151781702961021312, \\
&A^{\ast}_{639}=-641300074655810350440228116503844760152082975470094336\\
668 &a_{
55}=-5517441214194105433126639180879086430197314398412320, \\
&A^{\ast}_{663}=-58794819569722519246947052588417628381702293987381790208\\
692 &a_{
57}=-487293852416875056298891299076258026790715308037825216, \\
&A^{\ast}_{687}=-5382691989893920825057387379962379214878336698452619518976\\
716 &a_{
59}-43063324612408724058482351736525640971505158539537019360, \\
&A^{\ast}_{711}=-492331411094310243598851518968481327116177133949573219211776\\
740 &a_{
61}=-3807765135191026890629796012630986964392545971286263818880, \\
&A_{738}=-44999991655224865359801459729638972186485570857624927993968640\\
764 &a_{
63}=-336870497941041344941326929740498035221524766544197404568448, \\
&A^{\ast}_{759}=-4110683573025245915720418026614327134878535363601000658016970752
\\
788 &a_{
65}=-29817536107627026280375354657573625324749855885484211352178784, \\
&A^{\ast}_{783}=-3753105566589423622151865546742863862702362988456930306532547
79392\\
812 &a_{
67}=-2640483694133768694348825974762570078464246764217287044655040608, \\
 &A^{\ast}_{807}=-
34250066604647926889727509985278192843983793420163130299926828831232
\\
836 &a_{
69}=-2339302231720159785381190900580733565139572670910789181594822214
40, \\&A^{\ast}_{831}=-
3124209726760324826600719133662333295825193750667069283349549302802432\\
860 &a_{
71}=-2073334842464722080835335880678394765613836905527176962866019471
7760, \\&A^{\ast}_{855}=-
2848648403707089438800767077989207700176370570348417054372052052663531
52\\
884 &a_{
73}=-1838327467272503018023284563387674086088891665810123739273842888
409728, \\&A^{\ast}_{879}=-
2596376932291914595954593148339762930889889218217097595806788019761923
4816\\
\end{array}
\end{align*}

\hspace{-12pt}(iii)\ $n\equiv 20\pmod {24}$
\begin{align*}
\begin{array}{l|l}
n& a_{\alpha}, 
A^{\ast}_{[n/3]} \\ \hline 
908 &a_{
75}=-1630562611270581820339644567805649291670155784633977828245410059
09926944, \\&A^{\ast}_{903}=-
2365571801457082382449957806201062403238269316539050544966611715336538
277376\\
932 &a_{
77}=-1446787766179771143880641146327970387304955586054292653784397992
2706297728, \\
&A^{\ast}_{927}=-2154529731708362703498768412406995958086781390140282753533197
35457439571961856\\
956 &a_{
79}=-1284154078577872100670394856073262391516344033350714112670746082
005475905920, \\
&A^{\ast}_{951}=-1961663875672066714638183925465054355912396567314310941902880
9124847689381431296\\
\end{array}
\end{align*}

\end{appendix}


\begin{thebibliography}{MMM}

\bibitem[1]{1}
J. H. Conway, N.J.A. Sloane, ``Sphere Packings, Lattices and Groups," third edition, Springer, New York, 1999. 

\bibitem[2]{2}
C. L. Mallows, N. J. A. Sloane, An upper bound for self-dual codes, Inform. and Control 22 (1973) 188--200. 

\bibitem[3]{3}
S. Y. Zhang, On the nonexistence of extremal self-dual codes, Discrete Appl. Math. 91 (1999), 277--286. 

\end{thebibliography}
\end{document}